\documentclass[12pt]{article}

\usepackage{amsmath,amssymb}

\begin{document}

\pagestyle{myheadings} \markright{A CONJECTURAL LEFSCHETZ FORMULA...}

\title{A conjectural Lefschetz formula for locally symmetric spaces}
\author{Anton Deitmar\footnote{University of Exeter, 
Mathematics, Exeter EX4
4QE, England; a.h.j.deitmar@ex.ac.uk}}

\date{}
\maketitle

$$ $$

\def \1{{\bf 1}}
\def \a{{{\mathfrak a}}}
\def \ad{\mathop{\rm ad}\nolimits}
\def \al{\alpha}
\def \ar{{\alpha_r}}
\def \A{{\mathbb A}}
\def \Ad{\mathop{\rm Ad}\nolimits}
\def \Aut{{\rm Aut}}
\def \b{{{\mathfrak b}}}
\def \bs{\backslash}
\def \B{{\cal B}}
\def \c{{\mathfrak c}}
\def \cent{\mathop{\rm cent}\nolimits}
\def \C{{\mathbb C}}
\def \CA{{\cal A}}
\def \CB{{\cal B}}
\def \CC{{\cal C}}
\def \CD{{\cal D}}
\def \CE{{\cal E}}
\def \CF{{\cal F}}
\def \CG{{\cal G}}
\def \CH{{\cal H}}
\def \CHC{{\cal HC}}
\def \CL{{\cal L}}
\def \CM{{\cal M}}
\def \CN{{\cal N}}
\def \CP{{\cal P}}
\def \CQ{{\cal Q}}
\def \CO{{\cal O}}
\def \CS{{\cal S}}
\def \CT{{\cal T}}
\def \CV{{\cal V}}
\def \d{{\mathfrak d}}
\def \det{\mathop{\rm det}\nolimits}
\def \df{\ \begin{array}{c} _{\rm def}\\ ^{\displaystyle =}\end{array}\ }
\def \diag{\mathop{\rm diag}\nolimits}
\def \dist{\mathop{\rm dist}\nolimits}
\def \End{\mathop{\rm End}\nolimits}
\def \eps{\varepsilon}
\def \eqn{\begin{eqnarray*}}
\def \endeqn{\end{eqnarray*}}
\def \F{{\mathbb F}}
\def \Fx{{\mathfrak x}}
\def \FX{{\mathfrak X}}
\def \fin{_{\mathrm{fin}}}
\def \g{{{\mathfrak g}}}
\def \ga{\gamma}
\def \Ga{\Gamma}
\def \Gal{{\rm Gal}}
\def \GL{\mathop{\rm GL}\nolimits}
\def \h{{{\mathfrak h}}}
\def \Hom{\mathop{\rm Hom}\nolimits}
\def \im{\mathop{\rm im}\nolimits}
\def \Im{\mathop{\rm Im}\nolimits}
\def \Ind{\mathop{\rm Ind}\nolimits}
\def \k{{{\mathfrak k}}}
\def \K{{\cal K}}
\def\Kfin{K_{\rm fin}}
\def \l{{\mathfrak l}}
\def \la{\lambda}
\def \li{{\rm li}}
\def \La{\Lambda}
\def \m{{{\mathfrak m}}}
\def \n{{{\mathfrak n}}}
\def \name{\bf}
\def \Mat{\mathop{\rm Mat}\nolimits}
\def \N{\mathbb N}
\def \o{{\mathfrak o}}
\def \ord{\mathop{\rm ord}\nolimits}
\def \O{{\cal O}}
\def \p{{{\mathfrak p}}}
\def \ph{\varphi}
\def \prf{\noindent{\bf Proof: }}
\def \Per{{\rm Per}}
\def \PSL{{\rm PSL}}
\def \q{{\mathfrak q}}
\def \qed{\ifmmode\eqno $\square$\else\noproof\vskip 12pt plus 3pt minus 9pt \fi}
 \def\noproof{{\unskip\nobreak\hfill\penalty50\hskip2em\hbox{}%
     \nobreak\hfill $\square$\parfillskip=0pt%
     \finalhyphendemerits=0\par}}
\def \Q{\mathbb Q}
\def \res{\mathop{\rm res}\nolimits}
\def \R{{\mathbb R}}
\def \Re{\mathop{\rm Re}\nolimits}
\def \r{{\mathfrak r}}
\def \ra{\rightarrow}
\def \rank{\mathop{\rm rank}\nolimits}
\def \st{\mathop{\rm st}\nolimits}
\def \supp{\mathop{\rm supp}\nolimits}
\def \SL{\mathop{\rm SL}\nolimits}
\def \SO{\mathop{\rm SO}\nolimits}
\def \Spin{\mathop{\rm Spin}\nolimits}
\def \t{{{\mathfrak t}}}
\def \T{{\mathbb T}}
\def \tr{\mathop{\rm tr}\nolimits}
\def \vol{\mathop{\rm vol}\nolimits}
\def \z{\zeta}
\def \Z{\mathbb Z}
\def \={\ =\ }

\newcommand{\frack}[2]{\genfrac{}{}{0pt}{}{#1}{#2}}
\newcommand{\rez}[1]{\frac{1}{#1}}
\newcommand{\der}[1]{\frac{\partial}{\partial #1}}
\renewcommand{\binom}[2]{\left( \begin{array}{c}#1\\#2\end{array}\right)}
\newcommand{\norm}[1]{\left\|#1\right\|}
\renewcommand{\matrix}[4]{\left(\begin{array}{cc}#1 & #2 \\ #3 & #4 \end{array}\right)}
\renewcommand{\sp}[2]{\langle #1,#2\rangle}
\renewcommand{\labelenumi}{(\alph{enumi})}

\newtheorem{theorem}{Theorem}[section]
\newtheorem{conjecture}[theorem]{Conjecture}
\newtheorem{lemma}[theorem]{Lemma}
\newtheorem{corollary}[theorem]{Corollary}
\newtheorem{proposition}[theorem]{Proposition}

{ \tableofcontents}

\newpage
\section*{Introduction}
The theory of the Selberg Zeta Function is a vital part of the
spectral geometry of locally symmetric spaces, see \cite{Beilinson-Manin,
Bunke-Olbrich1, Bunke-Olbrich2, Cartier-Voros, Deitmar1, hitors, prod, det,
geom, Efrat, Efrat2, Elstrodt, Gangolli, Gangolli-Warner, Hejhal1, Hejhal2,
Juhl}.
For higher rank spaces there is no zeta function. If the space is compact, the
Lefschetz formula
\cite{hr, Juhl} seems to be an apropriate replacement as applications show
\cite{prime}. 

In this paper we suggest a Lefschetz formula for
non-compact finite volume spaces and we prove it in the case of Riemann
surfaces by exploiting the properties of the Selberg zeta function. This way
of proof might be extended to rank one spaces, but for higher rank a new idea
is required.

\section{Global Lefschetz numbers}
Let $G$ denote a connected semisimple Lie group with finite center. Fix a
maximal compact subgroup $K$ with Cartan involution $\theta$. So $\theta$ is
an automorphism of $G$ with $\theta^2=\rm Id$ and $K$ is the set of all $x\in
G$ with $\theta(x)=x$.

Let $\g_\R, \k_\R$ denote the real Lie algebras of $G$ and $K$ and let $\g$
and $\k$ denote their complexifications. This will be a general rule: for a
Lie group $H$ we denote by $\h_\R$ the Lie algebra of $H$ and by
$\h=\h_\R\otimes\C$ its complexification.
Let
$b:\g\times\g\to\C$ be a positive multiple of the Killing form. On $G,K$ and
all parabolic subgroups as well as all Levi-components we install Haar measures
given by the form $b$ as in
\cite{HC-HA1}.

Let $H$ be a non-compact Cartan subgroup of $G$. Modulo conjugation we can
assume that
$H=AB$ where $A$ is a connected split torus and $B$ is a closed subgroup of
$K$. Fix a parabolic $P$ with split component $A$. Then $P$ has Langlands
decomposition $P=MAN$ and $B$ is a Cartan subgroup of $M$. Note that  an
arbitrary parabolic subgroup $P'=M'A'N'$ of $G$ occurs in this way if and
only if the group $M'$ has a compact Cartan subgroup. In this case we say
that $P'$ is a \emph{cuspidal parabolic}.

The choice of the parabolic $P$ amounts to the same as a choice of a set of
positive roots $\Phi^+(\g,\a)$ in the root system $\Phi(\g,\a)$. The Lie
algebra $\n$ of the unipotent radical $N$ can be described as
$\n=\bigoplus_{\al\in\Phi^+(\g,\a)}\g_\al$, where $\g_\al$ is the \emph{root
space} attached to $\al$, i.e., $\g_\al$ is the space of all $X\in\g$ such
that $\ad(Y)X=\al(Y)X$ holds for every $Y\in\a$. Define
$\bar\n=\bigoplus_{\al\in\Phi^+(\g,\a)}\g_{-\al}$. This is the 
\emph{opposite} Lie algebra. Let $\bar\n_\R=\bar\n\cap\g_\R$ and $\bar
N=\exp(\bar\n_\R)$. Then $\bar P= MA\bar N$ is the \emph{opposite parabolic}
to $P$.

Let $\a^*$ denote the dual space of $\a$. Since $A=\exp(\a_0)$, every 
$\la\in\a^*$ induces a continuous homomorphism from $A$ to $\C^*$ written
$a\mapsto a^\la$ and given by $(\exp(H))^\la=e^{\la(H)}$.
Let $\rho_P\in\a^*$ be the half of the sum of all positive roots, each
weighted with its multiplicity. So $a^{2\rho_p}=\det(a|\n)$. Let
$\a_0^-\subset\a_0$ be the negative Weyl chamber consisting of all 
$X\in \a_0$ such that $\al(X)<0$ for every $\al\in\Phi^+(\g,\a)$. Let
$A^-=\exp(\a_0^-)$ be the negative Weyl chamber in $A$. Further let
$\overline{A^-}$ be the closure of $A^-$ in $A$. This is a manifold with
corners. Let $K_M=M\cap K$. Then $K_M$ is a maximal compact subgroup of $M$
and it contains $B$. Fix an irreducible unitary representation
$(\tau,V_\tau)$ of $K_M$. Then $V_\tau$ is finite dimensional. Let
$\breve\tau$ be the dual representation to $\tau$.

Let $\hat G$ denote the unitary dual of $G$, i.e., it is the set of all
isomorphy classes of irreducible unitary representations of $G$. Let $\hat
G_{\rm adm}\supset\hat G$ be the admissible dual. For
$\pi\in\hat G_{\rm adm}$ let $\pi_K$ denote the $(\g,K)$-module of $K$-finite
vectors in $\pi$ and let $\La_\pi\in\h^*$ be a representative of the
infinitesmal character of $\pi$. Let
$H^\bullet(\n,\pi_K)$ be the Lie algebra cohomology with coefficients in
$\pi_K$. By \cite{HeSch} for each
$q$ the
$(\a\oplus\m,K_M)$-module $H^q(\n,\pi_K)$ is admissible of finite length,
i.e., a Harish-Chandra module.

For $\la\in\a^*$ and an $A$-module $W$ let $W^\la$ denote the generalized
$\la$-eigenspace, i.e., $W^\la$ is the set of all $w\in W$ such that there is
$n\in\N$ with 
$$
(a-a^\la)^n w\= 0
$$
for every $a\in A$. 
Let $\m =\k_M\oplus\p_M$ be the Cartan decomposition of the Lie algebra $\m$
of $M$. For $\pi\in\hat G$ and $\la\in\a^*$ let $L_\la^\tau(\pi)$ denote the
representation-theoretic Lefschetz number given by
$$
L_\la^\tau(\pi)\df \sum_{p,q\ge 0}(-1)^{p+q+\dim
N}\dim\left(H^q(\n,\pi_K)^\la\otimes\wedge^p\p_M\otimes 
\breve\tau\right)^{K_M}.
$$

For a given smooth and compactly supported function $f\in
C_c^\infty(G)$ we define its Fourier transform $\hat f\colon \hat G\to\C$ by
$$
\hat f(\pi)\df \tr\pi(f).
$$

\begin{proposition}
\begin{enumerate}
\item
For every $\ph\in C_c^\infty(A^-)$ there exists $f_\ph\in C_c^\infty(G)$ such
that for every $\pi\in\hat G$,
$$
\widehat{f_\ph}(\pi)\=\sum_{\la\in\a^*} L_\la^\tau(\pi)\,\hat\ph(\la),
$$
where $\hat\ph$ is the Fourier transform of $\ph$, i.i.,
$\hat\ph(\la)=\int_A\ph(a)a^\la\,da$.
\item
The sum in (a) is finite, more precisely,
the Lefschetz number $L_\la^\tau(\pi)$ is zero unless there is an element $w$
of the Weyl group of $(\g,\h)$ such that
$$
\la\= \left.\left( w\La_\pi\right)\right|_\a -\rho_P.
$$
\end{enumerate}
\end{proposition}

\prf
The proof of part (a) is contained in section 4 of \cite{hr}, and
(b) is a consequence of Corollary 3.32 in \cite{HeSch}.
\qed

\section{Local Lefschetz numbers}
Let $\Ga\subset G$ be a discrete subgroup of finite covolume.
Let $X= G/K$ be the symmetric space and $X_\Ga=\Ga\bs X=\Ga\bs G/K$ be the
corresponding locally symmetric quotient.
The group $\Ga$ is called \emph{neat} if it is trosion-free and for every
finite dimensional representation $\rho$ of $G$ and every $\ga\in\Ga$ the
linear map $\rho(\ga)$ does not have a root of unity other than $1$ for an
eigenvalue. Every arithmetic group has a finite index subgroup which is
arithmetic and neat \cite{borel}. 

Let $L$ be a unimodular Lie group and $\Ga$ a lattice if $L$. Let $H\subset L$
be a Lie subgroup such that the Weyl-group $W(L,H)$ is finite. Here $W(L,H)$
is the normalizer of
$H$ modulo the centralizer of $H$. Let
$b$ be a non-degenerate symmetric bilinear form on the Lie algebra of $L$
which is invariant under
$H$. Suppose that there is a preferred Haar-measure $\mu_H$ on $H$.  
The form $b$ induces an $L$-invariant pseudo-Riemannian structire on $L/H$.
The Gau\ss -Bonnet construction (\cite{Dieu}, sect 24) extends to
pseudo-Riemannian structures to give an Euler-Poincar\'e measure $\eta$ on
$L/H$. Define a (signed) Haar-measure on $L$ by
$$
\mu_{b,H}\=\eta\otimes\mu_H.
$$
Define the $H$-index by
$$
\Ind_H(\Ga\bs L)\df \frac 1{|W(L,H)|}\mu_{b,H}(\Ga\bs L).
$$

\noindent
{\bf Remarks}
\begin{itemize}
\item
If $L$ is reductive, $H$ a compact Cartan subgroup and $\Ga$ cocompact and
torsion-free, then the $H$-indext equals the Euler-characteristic,
$$
\Ind_H(\Ga\bs L)\=\chi(\Ga\bs L/K_L),
$$ 
where $K_L$ is a maximal compact subgroup of $L$. 
\item
Assume $L$ reductive, $\Ga$ neat and $H=AB$ a Cartan subgroup with $A$
central in $L$. Let $C$ be the center of $L$, then $C=AB_C$, where $B_C$ is
compact. Let $\Ga_C=\Ga\cap C$ and $\Ga_A=A\cap \Ga_CB_C$ the projection of
$\Ga_C$ to $A$. Then $\Ga_A$ is a discrete and cocompact subgroup of $A$.
Under these circumstances,
$$
\Ind_H(\Ga\bs L)\=\vol(A/\Ga_A)\ \chi(A\Ga\bs L/ K_L).
$$
\item
Let $G$ as before and let $H=AB$ be a non-compact Cartan subgroup of $G$.
Let $\Ga\subset G$ be neat and let $[\ga]\in\CE_P(\Ga)$. Assume that $\ga$
is regular, then
$$
\Ind_H(\Ga_\ga\bs G_\ga)\= \vol(\Ga_\ga\bs G_\ga).
$$
\end{itemize}

Let $K_\ga\subset G_\ga$ be a maximal compact
subgroup and let $X_\ga=\Ga\cap G_\ga\bs G_\ga /K_\ga$ the corresponding
modular subvariety of $X_\Ga$. Being semisimple, the element $\ga$ lies in a
Cartan subgroup $H_\ga=A_\ga B_\ga$, where $A_\ga$ is a split connected
torus and $B$ is compact. Hence $\ga=\tilde a_\ga \tilde b_\ga$. If we
assume that $\tilde a_\ga$ is a regular element of $A_\ga$, then $A_\ga$ is
uniquely determined by $\ga$. 

Back to the notation of the first section let $\CE_P(\Ga)$ denote the set of
all conjugacy classes $[\ga]$ in $\Ga$ such that $\ga$ is in $G$ conjugate
to an element $a_\ga b_\ga$ of $A^- B$. Then there is a conjugate $H_\ga$ of
$H$ such that $\ga\in H_\ga$

For $[\ga]\in\CE_P(\Ga)$ we define the \emph{local Lefschetz number} by
$$
L^\tau(\ga)\df \Ind_{H_\ga}(\Ga_\ga\bs
G_\ga)\frac{\tr\tau(b_\ga)}{\det(1-a_\ga b_\ga |\n)}.
$$

\section{The Lefschetz formula}
The unitary $G$-representation on $L^2(\Ga\bs G)$ decomposes as
$$
L^2(\Ga\bs G)\= L_{\rm disc}^2\oplus L_{\rm cont}^2,
$$
where 
$$
L_{\rm disc}^2\= \bigoplus_{\pi\in\hat G}N_\Ga(\pi)\,\pi
$$
is a direct sum of irreducibles with finite multiplicities and $L_{\rm
cont}^2$ is a sum of continuous Hilbert integrals. In particular, $L_{\rm
cont}^2$ does not contain any irreducible subrepresentation.

Let $r$ be the dimension of $A$ and let $\al_1,\dots,\al_r\in \a_\R^*$ be
the primitive positive roots. Let $\a_\R^{*,+}=\{ t_1\al_1+\cdots +t_r\al_r
: t_1,\dots,t_r>0\}$ be the positive dual cone and let
$\overline{\a_\R^{*,+}}$ be its closure in $\a_\R^*$.

For $\mu\in\a^*$ and $j\in\N$ let $C^{\mu,j}(A^-)$ denote the space of all
functions on $A$ which
\begin{itemize}
\item
are $j$-times continuously differentiable on $A$,
\item
are zero outside $A^-$,
\item
satisfy $|D\ph|\le C |a^\mu|$ for every invariant diffferential operator $D$
on $A$ of degree $\le j$, where $C>0$ is a constant, which 
depends on $D$. 
\end{itemize}

This space can be topologized with the seminorms
$$
N_D(\ph)\= \sup_{a\in A} |a^{-\mu}D\ph(a)|,
$$
$D\in U(\a)$, $\deg(D)\le j$. Since the space of operators $D$ as above is
finite dimensional, one can choose a basis $D_1,\dots,D_n$ and set
$$
\norm{\ph}\= N_{D_1}(\ph)+\cdots +N_{D_n}(\ph).
$$
The topology of $C^{\mu,j}(A^-)$ is given by this norm and thus
$C^{\mu,j}(A^-)$ is a Banach space.

\begin{conjecture}\label{Lefschetz}
(Lefschetz Formula)\\
For $\la\in \a*$ and $\pi\in\hat
G_{\rm adm}$ there is an integer $N_{\Ga,\rm cont}(\pi,\la)$ which
vanishes if $\Re(
\la)\notin \overline{\a_{\R}^{*,+}}$ and there are $\mu\in\A^*$ and
$j\in\N$ such that for each
$\ph\in C^{\mu,j}(A^-)$ and with
$$
m_\la(\pi)\df N_\Ga(\pi)+N_{\Ga,\rm cont}(\pi,\la)
$$
we have
$$
\sum_{\stackrel{\pi\in\hat G}{\la\in\a^*}} 
m_\la(\pi)\,
L_\la^\tau(\pi)\,\int_A\ph(a) a^\la\, da
\=
 \sum_{[\ga]\in\CE_P(\Ga)} L^\tau(\ga)\,\ph(a_\ga).
$$
Either side of this identity represents a continuous functional on
$C^{\mu,j}(A^-)$.
\end{conjecture}

In the following cases the conjecture is known.
\begin{enumerate}
\item
The conjecture holds if $\Ga$ is cocompact.
In that case the numbers $N_{\Ga,\rm cont}(\pi,\la)$ are all zero. This is
shown in \cite{hr}.
\item
In the next section we will prove the conjecture for $G=\PSL_2(\R)$.
\end{enumerate}

We will now make the conjecture more precise for congruence subgroups. For
this assume that $G=\CG(\R)$ for some semisimple linear group $\CG$ defined
over $\Q$. Let $\A=\A_{\rm fin}\times\R$ be the adele ring over $\Q$.
Assume that $\Ga$ is a congruence subgroup., i.e., there
exists a compact open subgroup $K_\Ga$ of $\CG(\A_{\rm
fin})$ such that $\Ga=\CG(\Q)\cap K_\Ga$. To explain the conjectured nature
of the number $N_{\Ga \rm cont}(\pi)$ we will recall Arthur's trace
formula. This formula is the equality of two distributions on $\CG(\A)$,
$$
J_{\rm geom}\= J_{\rm spec}.
$$
The geometric distribution $J_{geom}$ can be described in terms of weighted
orbital integrals. Our interest however is focused on the spectral
distribution $J_{\rm spec}$.
According to~\cite{Art-eisII}, Theorem~8.2, one has
$$
J_{\mathrm{spec}}(f)=\sum_\chi J_\chi(f),
$$
where $\chi$ runs through conjugacy classes of pairs
$(\CM_0, \pi_0)$ consisting of a $\Q$-rational Levi
subgroup $\CM_0$ and its cuspidal automorphic
representation~$\pi_0$, the sum being absolutely convergent.
The particular terms have expansions
$$
J_\chi(f)=\sum_{\CM,\eta} J_{\chi,\CM,\eta}(f),
$$
where the sum runs over all $\Q$-rational Levi
subgroups $\CM$ of $\CG$ containing a fixed minimal one
(which we take to be the subgroup $\CA_0$ of diagonal matrices)
and, for each~$\CM$, over all automorphic representations $\eta$ of
$\CM(\A)^1$. Explicitly,
$$
J_{\chi,\CM,\eta}(f)=\sum_{s\in W_\CM}
c_{\CM,s}\int_{i(\a_{\CL}^{\CG})^*}\sum_{\CP}
\tr\left(\mathfrak M_{\CL}(\CP,\nu)M(\CP,s)
    \rho_{\chi,\eta}(\CP,\nu,f)\right)d\nu.
$$
Here, for a given element $s$ of the Weyl group of $\CM$ in $\CG$,
the Levi subgroup $\CL$ is determined by $\a_\CL=(\a_\CM)^s$,
and $\CP$ runs through all parabolic subgroups of $\CG$ having $\CM$
as a Levi component. 
Let us comment on the items in the integrand. Let
$\rho(\CP,\nu)$ be the representation of $\CG(\A)$ which
is induced from the representation
of $\CP(\A)$ in
$$
L^2(\CM(\Q)\bs \CM(\A))\cong L^2(\CN(\A)\CP(\Q)\bs \CP(\A))
$$
twisted by~$\nu$. If one starts the induction with the
subspace of the $\pi$-isotypical component spanned by
certain residues of Eisenstein series coming from $\chi$,
one gets a subrepresentation which is denoted by
$\rho_{\chi,\eta}(\CP,\nu)$. We let $\rho_{\chi,\eta}(\CP,\nu,f)$
act in the space of $\rho(\CP,\nu)$ by composing it with the
appropriate projector. Further, there is a
meromorphic family of standard intertwining operators
$M_{\CQ|\CP}(\nu)$ between dense subspaces of
$\rho(\CP,\nu)$ and $\rho(\CQ,\nu)$ defined by an integral
for $\Re\nu$ in a certain chamber. The operator $M(\CP,s)$
is $M_{s\CP|\CP}(0)$ followed by translation with a
representative of $s$ in $\CG(\Q)$. And finally,
$\mathfrak M_{\CL}(\CP,\nu)$ is obtained from such intertwining
operators by a limiting process.
The operator valued function 
$$
\nu\ \mapsto\ \frak M_\CL(\CP,\nu)M(\CP,s)
$$
extends to a meromorphic function on $(\a_\CL^\CG)^*$. For $\nu\in
(\a_\CL^\CG)^*$ let $R_\nu$ denote the residue of this operator valued
function at $\nu$. Arthur proved that the distribution
$$
f\ \mapsto\ \tr\left( R_\nu\,\rho_{\chi,\eta}(\CP,\nu,f)\right) = D(f)
$$
is invariant. Let $\1_{K_\Ga}$ be the indicator function of $K_\Ga$. We
conjecture that the distribution on $G$,
$$
D_\infty \colon\ \ph\ \mapsto\ D\left( \frac
1{\vol(K_\Ga)}\1_{K_\Ga}\otimes\ph\right)
$$
is a finite linear combination of traces with integer coefficients, i.e.,
$$
D_\infty(\ph)\=\sum_{\eta\in\hat G_{\rm adm}}
c(\chi,\eta,\CP,\nu,\pi)\,\tr\pi(\ph)
$$
for some $c(\chi,\eta,\CP,\nu,\pi)\in\Z$.

\begin{conjecture}
Conjecture \ref{Lefschetz} holds with
$$
N_{\Ga,\rm cont}(\pi,\la)\=
\1_{\overline{\a_\R^{*,+}}}(\la)\sum_{\chi,\eta,\CP,\nu}
c(\chi,\eta,\CP,\nu,\pi).
$$
\end{conjecture}

\section{$\PSL_2(\R)$}
In this section we will prove the conjecture in the simplest case, that of
the group $G=\PSL_2(\R)=\SL_2(\R)/\pm 1$. For this group there is, up to
conjugation, only one choice for
$H=AB$, namely $A=\left\{ \pm\matrix t{}{}{1/t}: t>0\right\}$ and
$B=M=\{1\}$. We further choose $N=\left\{\pm\matrix 1x{}1 :x\in\R\right\}$
and set
$P=MAN=AN$. Fix the maximal compact subgroup $K={\SO}_2(\R)/\pm 1$. Finally,
we choose the from
$b$ to be
$b(X,Y)=\tr(XY)$. Occasionally we will view a function $f$ on $G$ as
a function on $\SL_2(\R)$ with $f(-x)=f(x)$.

Let us recall some facts from the representation theory of
$\PSL_2(\R)$. Let $\la\in\a^*$
and denote by $\pi_{\la}$ the corresponding principal series
representation, so $\pi_{\la}$ lives on the space of measurable
functions $f\colon G\to\C$ with $f(an x)= a^{\la+\rho} f(x)$ which
are square integrable on $K$, modulo nullfunctions. The representation is
the right regular representation, i.e., $\pi(y)f(x)=f(xy)$. For each natural
number
$n$ there are exact sequences
$$
0\ \to\ \CD_{2n}^+\oplus\CD_{2n}^-\ \to\ \pi_{(2n-1)\rho}\
\to\ \delta_{2n-1}\ \to\ 0,
$$
and
$$
0\ \to\  \delta_{2n-1}\ \to\ \pi_{(1-2n)\rho}\ \to\
\CD_{2n}^+\oplus\CD_{2n}^-\ \to\ 0,
$$
where $\CD_{2n}^{\pm}$ are the discrete series, resp. limit of discrete series
representations as in \cite{Knapp}, and $\delta_{2n-1}$ is the unique
irreducible representation of $G$ of dimension $2n-1$.
In all other cases the representation $\pi_{\la}$ is
irreducible. If  $\la$ is purely imaginary, then $\pi_{\la}$ is
isomorphic with
$\pi_{-\la}$ and this is the only isomorphism between different
principal series representations.
The admissible dual $\hat G_{\rm adm}$ consists of all irreducible principal
series representations and all $\CD_{2n}^\pm$ and all $\delta_{2n-1}$.
The unitary dual $\hat G$ comprises all irreducible $\pi_{\la}$, where
$\la$  is purely imaginary, all $\pi_{t\rho}$ for $0<t<1$ and all
$\CD_{2n}^\pm$. For $\la\in\a^*$  we also write
$\la$ for the quasi-character $a\mapsto a^\la$ of the
group $A$.

\begin{proposition}
\begin{enumerate}
\item
If $\pi\in\hat G_{\rm adm}$ is a principal series representation, then
$H^0(\n,\pi_K)=0$. More generally, we have $H^0(\n,(\pi_{\la})_K)=0$
unless $\la=(1-2k)\rho$ for some $k\in\N$ in which case it is
one-dimensional and $A$ acts via $(2k-2)\rho$.
\item
$H^0(\n,\delta_{2n-1})$ is one dimmensional for every $n\in\N$ and $A$ acts
via the character
$(1-2n)\rho$.
\item
$H^0(\n,\CD_{2n}^\pm)=0$ for every $n\in\N$.
\end{enumerate}
\end{proposition}

\prf
Recall the Iwasawa decomposition $G=ANK$.\\ Explicitly, for $g=\pm\matrix
abcd\in G$ we get
$$
g\=\pm\matrix{\frac 1{\sqrt{c^2+d^2}}}{}{}{\sqrt{c^2+d^2}}
\matrix 1{ac+bd}{}1
\frac 1{\sqrt{c^2+d^2}}\matrix d{-c}cd.
$$
Let $f\in H^0(\n,\pi_K)=\pi_K^\n$. Then for every $n\in N$ we have
$f(xn)=f(x)$. Let $w=\pm\matrix {}{-1}1{}$. If $f(w)=0$ then $f=0$. So we
assume that $f(w)=1$. Let $\la=s\rho$, $s\in\C$. Then for $x\in\R$,
\begin{eqnarray*}
1 &=& f\left(w\matrix 1x{}1\right)\\
&=& f\matrix 0{-1}1x\\
&=& f\left(\matrix{\frac 1{\sqrt{1+x^2}}}{}{}{\sqrt{1+x^2}}\matrix
1{-x}{}1 \frac 1{\sqrt{1+x^2}}\matrix x{-1}1x\right)\\
&=& \sqrt{1+x^2}^{-s-1}
f\left(\frac 1{\sqrt{1+x^2}}\matrix x{-1}1x\right),
\end{eqnarray*}
or
$$
f\left(\pm\frac 1{\sqrt{1+x^2}}\matrix x{-1}1x\right)
\=\sqrt{1+x^2}^{s+1}.
$$
Since $f\in\pi_K$, the function $f$ is continuous on $K$, so the limit as
$x\to\infty$ must exist, which implies $s=-1$ or $\Re(s)<-1$.
In the case $s=-1$ the constant function $f(x)=1$ indeed gives a basis for
$H^0(\n,\pi_K)$. If $\Re(s)<-1$ then $f$ can only be a smooth function on
$K$ if $s=1-2k$ is integral and odd. In order to determine the $A$-actions we
need to introduce some notation. For a $\C[A]$-module $V$ and
$\la\in\a^*$ we write $V^\la$ for the generalized $\la$-eigenspace in $V$.
This means $v\in V^\la$ if and only if there is $n\in\N$ such that
$$
(a-a^\la)^n v\= 0
$$
for every $a\in A$. On $\a^*$ we introduce a partial order $<$ as follows,
$$
\nu <\mu\ \Leftrightarrow\ \mu-\nu\in\N\rho.
$$

\begin{lemma}
\begin{enumerate}
\item
For $p=0,1$ and $\pi\in\hat G_{\rm adm}$ we have
$$
H^p(\n,\pi_K)\=\bigoplus_{\nu=w\Lambda_\pi|_\a}H^p(\n,\pi_K)^{\nu-\rho},
$$
where $\La_\pi\in\h^*$ is a representative of the infinitesimal character of
$\pi$ and the sum runs over $w$ in the Weyl group $W(\g,\h)$. 
\item
If $H^0(\n,\pi_K)^\mu\ne 0$, then there is $\nu\in\a^*$ with $\nu <\mu$ such
that $H^1(\n,\pi_K)^\nu\ne 0$.
\end{enumerate}
\end{lemma}

\prf
Part (a) is Corollary 3.32 of \cite{HeSch} and part (b) is Proposition 2.32
of \cite{HeSch}.
\qed

Let $\pi=\pi_{(1-2k)\rho}$. According to part (a) of the Lemma, the group
$A$ acts on $H^0(\n,\pi_K)$ either via $(2k-2)\rho$ or $-2k\rho$. The
second possibility is excluded by part (b) of the Lemma. This proves part
(a) of the Proposition.

Part (b) of the Proposition follows from highest weight theory and part (c)
follows from part (a) and the fact that $\CD_{2n}^+\oplus\CD_{2n}^-$ is a
subrepresentation of
$\pi_{2n-1}$.
\qed

\begin{proposition}
For $\pi\in\hat G_{\rm adm}$ the dimension of $H^1(\n,\pi_K)=1$ is one except
if
$\pi=\pi_{\la}$ where $\la$ is purely imaginary (unitary principal
series) in which case the dimension is two.

If $\pi=\pi_{\la}$ is nonunitary principal series, then $A$ acts via
$\la-\rho$.

If $\pi=\pi_{\la}$ is unitary principal series, then A acts via
$(\la-\rho)\ \oplus\ (-\la-\rho)$.

If $\pi=\delta_{2n-1}$, then A acts via $-2n\rho$.

If $\pi=\CD_{2n}^\pm$, then A acts via $(2n-2)\rho$.
\end{proposition}

\prf
For any $MA$-module $U$ we have
$$
\Hom_{MA}(H^1(\n,\pi_K),U\otimes\C_{-\rho})\=\Hom_G(\pi,\Ind_P^G(U)),
$$
(see Theorem 4.9 of \cite{HeSch}). The Proposition follows from this.
\qed

Since there is no choice for $\tau$ we leave this index out of the notation
for the Lefschetz numbers.

\begin{proposition}
\begin{enumerate}
\item
Let $\pi_{\mu}\in\hat G_{\rm adm}$ be a non-unitary principal series
representation. Then $L_\la(\pi_{\la})=0$ unless 
$\la=\mu-\rho$. In that case,
$$
L_{\mu-\rho}(\pi_{\mu})\= 1.
$$
\item
Let $\pi_{\la}$ be a unitary principal series. Then
$$
L_{\mu-\rho}(\pi_{\mu})\= 1\= L_{-\mu-\rho}(\pi_{\mu}),
$$
and $L_\la(\pi_{\mu})=0$ in all other cases.
\item
Let $n\in\N$. Then $L_\la(\delta_{2n-1})=0$ except for
$$
L_{(1-2n)\rho}(\delta_{2n-1})\= -1,
$$
and
$$
L_{-2n\rho}(\delta_{2n-1})\= 1.
$$
\item
Let $n\in\N$. Then $L_\la^\tau(\CD_{2n}^\pm)=0$ except for
$$
L_{(2n-2)\rho}(\CD_{2n}^\pm)\= 1.
$$
\end{enumerate}
\end{proposition}

Next for the local Lefschetz numbers. If $[\ga]\in\CE_P(\Ga)$, then $\ga$ is
$G$-conjugate to $\pm\matrix {N(\ga)^{1/2}}{}{}{N(\ga)^{-1/2}}$ for some
$N(\ga)>1$. An element $\ga$ of $\Ga$ is called \emph{primitive} if
$\ga=\sigma^n$ for
$\sigma\in\Ga$ and $n\in\N$ implies $n=1$. Each $\ga\in\CE_P(\Ga)$ is a
power of a unique primitive $\ga_0$ which will be called the primitive
\emph{underlying} $\ga$.

We write $L(\ga)$ for $L^\tau(\ga)$ since $\tau $ is trivial anyway. Then
$$
L(\ga)\=\frac{\log N(\ga_0)}{1-N(\ga)^{-1}}.
$$

We will now recall some facts about the \emph{Selberg zeta function}
\cite{Efrat, Hejhal, Iwaniec}. Let $\CE_P^p(\Ga)$ denote the set of all
primitive classes in $\CE_P(\Ga)$. The Selberg zeta function is given by
the product
$$
Z(s)\=\prod_{\ga\in\CE_P^p(\Ga)}\prod_{k=0}^\infty \left(
1-N(\ga)^{-s-k}\right).
$$
The product converges locally uniformly for $\Re(s)>1$. The zeta function
extends to a meromorphic function on the plane of finite order. It has a
simple zero at $s=1$ and zeros at $s=\frac 12 \pm u$ of multiplicity
$N_\Ga(\pi_{u\frac\rho 2})$. These are all zeros or poles in $\Re(s)\ge
\frac 12$ except for $s=\frac 12$ where $Z(s)$ has a zero or pole of order
$N_\Ga(\pi_0)$ minus the number of cusps. The poles and zeros in
$\Re(s)<\frac 12$ can be described through the scattering matrix or
intertwining operators \cite{Efrat, Hejhal, Iwaniec}.

Recall the inversion formula for the Mellin transform. Let the
function $\psi$ be integrable on $(0,\infty)$ with respect to the
measure $\frac{dt}t$, in other words, $\psi\in L^1\left((0,\infty),
\frac{dt}t\right)$. Then the \emph{Mellin transform} of $\psi$ is given by
$$
M\psi(s)\df \int_0^\infty t^s\,\psi(t)\,\frac{dt}t,\qquad s\in i\R.
$$
If $\psi$ is continuously differentiable and $\psi'(t)t, \psi''(t)t^2$ are
also in $L^1\left((0,\infty),\frac{dt}t\right)$, then the following
inversion formula holds,
$$
\psi(t)\= \frac{1}{2\pi i} \int_{i\R} M\psi(s) t^{-s}\, ds.
$$
Now assume that $\psi$ is supported in the interval $[1,\infty)$ and that
for some $\mu>0$ the functions $\psi(t),\psi'(t)t,\psi''(t)t^2$ all are
$O(t^{-\mu})$. Then it follows that the integral
$M\psi(s)$ defines a function holomorphic in
$\Re(s)<\mu$ and the integral in the
inversion formula can be shifted,
$$
\psi(t)\= \frac{1}{2\pi i} \int_{C-i\infty}^{C+i\infty} M\psi(s)\, t^{-s}\,
ds,
$$
for every $C<\mu$.

Every $\ga\in\CE_P(\Ga)$ can be
written as
$\ga=\ga_0^n$ for some uniquely determined $\ga_0\in\CE_P^p(\Ga)$ and a
unique $n\in\N$.
A computation yields for $\Re(s)>1$,
\begin{eqnarray*}
\frac {Z'}Z(s)&=& \sum_{\ga\in\CE_P^p(\Ga)}\sum_{n=1}^\infty \frac{\log
N(\ga)}{1-N(\ga)^{-n}}\, N(\ga)^{-ns}\\
&=& \sum_{\ga\in\CE_P(\Ga)} \frac{\log
N(\ga_0)}{1-N(\ga)^{-1}}\, N(\ga)^{-s}
\end{eqnarray*}
Let $\psi$ be as above with $\mu>1$ and let $1<C<\mu$. Then, since
$\frac{Z'}Z(s)$ is bounded in $\Re(s)=C$ we can interchange integration and
summation to get
\begin{eqnarray*}
\frac 1{2\pi i}\int_{C-i\infty}^{C+i\infty} \frac{Z'}Z(s)\, M\psi(s)\, ds
&=& \sum_{\ga\in\CE_P(\Ga)}\frac{\log N(\ga_0)}{1-N(\ga)^{-1}}\psi(N(\ga)).
\end{eqnarray*}\vspace{10pt}

For $a\in A^-=\left\{\matrix t{}{}{t^{-1}} : 0<t<1\right\}$ set
$$
\ph(a)\= \ph\matrix t{}{}{t^{-1}}\df \psi\left(\frac 1t\right).
$$
Then $\ph\in C^{2,2\mu\rho}(A^-)$ and 
\begin{eqnarray*}
\frac 1{2\pi i}\int_{C-i\infty}^{C+i\infty}\frac{Z'}Z(s)\, M\psi(s)\, ds &=&
\sum_{\ga\in\CE_P(\Ga)}\frac{\log N(\ga_0)}{1-N(\ga)^{-1}}\ph(a_\ga)\\
&=& \sum_{\ga\in\CE_P(\Ga)}L(\ga)\ph(a_\ga),
\end{eqnarray*}
which is the right hand side of the Lefschetz formula.

Now suppose that $\ph\in C^{j,2\mu\rho}(A^-)$ for some $j\in\N$ and some
$\mu>1$. Then the functions $\psi(t),\psi'(t)t,\dots,\psi^{(j)}(t)t^j$ are
all $O(t^{-\mu})$. Integration by parts shows that
$$
M\psi(t)\-\frac{(-1)^j}{s(s+1)\cdots (s+j-1)}\int_0^\infty
t^{s}\,\psi^{(j)}(t)t^j\,\frac {dt}t.
$$
This implies that $M\psi(s)=O\left((1+|s|)^{-j}\right)$ uniformly in
$\{\Re(s)\le\al\}$ for every $\al <\mu$.

For $R>0$ and $a\in\C$ let $B_r(a)$ be the closed disk around $a$ of radius
$r$. Let $g$ be a meromorphic function on $\C$ with poles $a_1, a_2,\dots$.
We say that $g$ is \emph{essentially of moderate growth}, if there is a
natural number $N$, a constant $C>0$, and  as sequence of real numbers
$r_n>0$ tending to zero, such that the disks $B_{r_n}(A_n)$ are pairwise
disjoint and that on the domain $D=\C\setminus\bigcup_n B_{r_n}(a_n)$ it
holds $|g(z)|\le C|z|^N$. Every such $N$ is called a 
\emph{growth exponent}
of $g$.

\begin{lemma}
Let $h$ be a meromorphic function on $\C$ of finite order and let $g=h'/h$
be its logarithmic derivative. Then $g$ is essentially of moderate growth
with growth exponent equal to the order of $h$ plus two.
\end{lemma}

\prf
This is a direct consequence of Hadamard's factorization Theorem applied to
$h$.
\qed

This Lemma together with the growth estimate for $M\psi$ implies that for
$j$ large enough the contour integral over $C+i\R$ can be moved to the left,
deforming it slightly, so that one stays in the domain $D$, and gathering
residues. Ultimately, the contour integral will tend to zero, leaving only
the residues. One gets
\begin{eqnarray*}
\sum_{\ga\in \CE_P(\Ga)} L(\ga)\, \ph(a_\ga) &=& \sum_{s_0\in\C}\left(
\res_{s=s_0}\frac{Z'}Z(s)\right)\, M\psi(s_0)\\
&=& \sum_{s_0\in\C}\left(
\res_{s=s_0}\frac{Z'}Z(s)\right)\, \int_0^\infty\psi(t) t^{s_0}\frac{dt}t\\
&=&\sum_{s_0\in\C}\left(
\res_{s=s_0}\frac{Z'}Z(s)\right)\, \int_{A^-}\ph(a) a^{-s_0\rho}\, da.
\end{eqnarray*}
This implies the conjecture in the case $G=\PSL_2(\R)$.

\section{Applications}
If we assume the Lefschetz formula in general, then most applications known in
the compact case carry over to the non-compact case. In this section we will
only highlight the prime geodesic theorem as an example.
For this we assume that the parabolic $P=MAN$ is \emph{minimal}, i.e., the
group $M$ is compact. Let $r=\dim A$ and let $\al_1,\dots,\al_r$ be positive
multiples of simple roots such that for the modular shift $\rho_P$ we have
$$
2\rho_P\=\al_1+\cdots +\al_r.
$$
This fixes $\al_1,\dots,\al_r$ up to order. We choose a Haar-measure (i.e., a
form $b$) such that for the subset of $A$,
$$
\{ a\in A : 0\le \al_k(\log a)\le 1,\  k=1,\dots,r\}
$$
has volume 1.

\begin{theorem} (Prime Geodesic Theorem)\\
For $T_1,\dots,T_r >0$ let
$$
\Psi(T_1,\dots,T_r)\=\sum_{\stackrel{\stackrel{[\ga]\in\CE_P(\Ga)}{}}{\stackrel{}{a_\ga^{-\al_k}\le
T_k,\ k=1,\dots,r}}}
\la_\ga.
$$
We assume that the Lefschetz formula holds.
Then, as $T_k\to\infty$ for $k=1,\dots,r$ we have
$$
\Psi(T_1,\dots,T_r)\ \sim\ T_1\cdots T_r.
$$
\end{theorem}

\prf The proof of the compact case \cite{prime} carries over.
\qed

We further note a consequence of the Prime Geodesic Theorem which comes about
when one applies the Prime geodesic Theorem to $G=\SL_d(\R)$ and
$\Ga=\SL_d(\Z)$. Let
$d$ be a prime number $\ge 3$. Let $\CC$ be the set of all totally real number
fields
$F$ of degree $d$. Let $O$ be the set of all orders $\CO$ in number fields
$F\in\CC$. For an order $\CO\in O$ let $h(\CO)$ be its class number and
$R(\CO)$ its regulator. For $\la\in\CO^\times$ let $\sigma_1,\dots,\sigma_d$
denote the real embeddings of $F$ ordere in a way that
$|\sigma_k(\la)|\ge|\sigma_{k+1}(\la)|$ holds for $k=1,\dots, d-1$. For $k$ in
the same range let
$$
\al_k(\la)\df k(d-k)\,\log\left(
\frac{|\rho_{k}(\la)|}{|\rho_{k+1}(\la)|}\right).
$$
Let
$$
c=(\sqrt 2)^{1-d}  \left(\prod_{k=1}^{d-1}2k(d-k)\right).
$$
 So $c>0$ and it comes about as
correctional factor between the Haar measure normalization used in the Prime Geodesic
Theorem and the normalization used in the definition of the regulator.

\begin{theorem}\label{A.1}
For $T_1,\dots,T_r>0$ set
$$
\vartheta(T)\df \sum_{\stackrel{\la\in\CO^\times/\pm
1,\ \CO\in O}{\stackrel{0<\al_k(\la)\le T_k}{k=1,\dots,d-1}}}  R(\CO)\,
h(\CO).
$$
Then we have, as $T_1,\dots,T_{d-1}\ra \infty$,
$$
\vartheta(T_1,\dots,T_{d-1})\ \sim\ \frac{c}{\sqrt{d}}\,T_1\cdots T_{d-1}.
$$
\end{theorem}

\end{document}